\documentclass[conference,12pt,draftcls,onecolumn]{IEEEtran}
\usepackage[english]{babel}
\usepackage[utf8]{inputenc}
\usepackage{amsmath}
\usepackage{graphicx}
\usepackage[usenames, dvipsnames]{color}
\usepackage{amssymb}
\usepackage{xfrac}
\usepackage{amssymb}
\usepackage{verbatim}
\usepackage{url}
\usepackage{epstopdf}
 \usepackage{lettrine}   

\newcommand*{\QEDB}{\hfill\ensuremath{\square}}%

\newcommand*{\doverline}[1]{\overline{\overline{#1}}}
\graphicspath{{Figures/}}
\newcommand{\matr}[1]{\mathbf{#1}}

\DeclareMathOperator*{\argmax}{arg\,max}
\DeclareMathOperator*{\argsup}{arg\,sup}
\ifCLASSINFOpdf
\else
\fi
\begin{document}
%

\title{Cooperation and Competition among Energy Storages}

\author{Jesus E. Contreras-Oca\~na,~\IEEEmembership{Student Member,~IEEE,}
        Miguel A. Ortega-Vazquez,~\IEEEmembership{Senior Member,~IEEE,}
        Baosen Zhang,~\IEEEmembership{Member,~IEEE}
 \thanks{Jesus~E.~Contreras~Oca\~na, Miguel~A.~Ortega-Vazquez and Baosen Zhang are with the Department
of Electrical Engineering at the University of Washington. Emails: \{ jcontrer, maov,zhangbao\}@uw.edu. This work is partially supported by the University of Washington Clean Energy Institute.}
}

%


\maketitle

\begin{abstract}
We study competition and cooperation among a group of storage units. As the number of energy storages increases, the profit of storages approaches zero under competition.  We propose two ways in which storages can achieve the maximum possible profit.  The first is a decentralized approach in which storages incur artificial costs that act as incentives for them to behave as a coalition. No private information needs to be exchanged between the storages to calculate the artificial cost function. The second is a centralized approach in which an aggregator coordinates and splits profits with storages in order to achieve maximum profit. We do not assume the nature of the storage-aggregator relationship and derive the necessary conditions for longterm cooperation. We use Nash's axiomatic bargaining problem to model and predict the profit split between aggregator and storages. 
\end{abstract}


%
\IEEEpeerreviewmaketitle



\section{Introduction}
\lettrine{L}{arge} scale introduction energy storage to the grid has the potential to increase the efficiency of the power system from several dimensions: by shifting load from low to high price hours, providing reserves, improving power quality, and even defering capital investments \cite{Denholm}.  Energy storage is of particular importance to the successful integration of renewable sources into the power system as it can be used to mitigate their well-known stochasticity and intermittency. Storage devices range from large pumped hydro plants to household-level storage units (\emph{e.g.} Tesla Powerwall \cite{Tesla}). Although the amount of energy storage in the grid is currently limited, it has recently undergone an unprecedented growth and is expected to continue doing so as costs are driven down \cite{IRENA}.  From a storage owner perspective, these devices are capable of arbitraging energy in time, by buying energy at low prices and selling it back later at higher prices.  

The optimal utilization of storage devices has been explored under both regulated and competitive market environments.  In the former, storage is seen as a public asset that is centrally operated by a system operator (SO) to minimize the operating cost of the system \cite{Abdurrahman,Kinter-Meyer,Ortega-Vazquez_TPS_May_2005,Pozo}.   The latter entails a decentralized operation where storage devices pursue their own objective (\emph{e.g.} maximize profits) in a market environment \cite{Sarker_M}. 

Storages have been treated as price-takers if their capacity is insufficient to re-shape the system demand \cite{Figueiredo,Walawalkar}.  However, as the number of storages in the system increases, its impact on the residual demand and electricity prices increases \cite{Jenkin,Sioshansi}. Therefore storage units would need to be modeled as \emph{price-anticipatory} units-- units act under the assumption that the other players' actions influence prices.  

Most of research have focused on eliminating market power that arise when storages are price-anticipatory~\cite{Sioshansi_EJ_2010}.   In \cite{He} the authors propose a framework in which storage owners auction physically binding rights to their storage capacity.  In \cite{Taylor}, it is proposed for the storage owners to sell financially binding rights via a market operated by the SO.  Finally, the authors of \cite{Munoz-Alvares} propose a framework in which the storage devices are treated as a communal asset. 

There are two common threads in these approaches. The first is that they require a third party operating the storage devices, which might not necessarily be in the best interests of the storage owner. The second, and arguably more important, is that these approaches tend to drive storages out of the market. By eliminating power power, storages tend to make zero profit as their number increases, and may lead them to leave the market all together. 

We take a different viewpoint by encouraging storages to form \emph{coalitions} and thus achieve maximum profitability. We adopt the maximum profit as the objective for two main reasons.  Firstly, the social welfare is not easily defined at times. For instance, the maximum social welfare for all the storage units in a utility service area depends on the objective function of the utility, which does not necessary represent the least cost solution. Secondly, maximizing the profit of storage, they are encouraged to stay in the system and may lead to faster adoption of storage technologies. 


This paper studies the arbitrage problem in a market setting, where multiple distributed storage units compete in a price-anticipatory manner. Due to competition, each unit's profits is lower than the profit they would obtain if they had cooperated. The contributions of this paper are two mechanisms, one distributed and one centralized, that can be used to incentivize cooperation among a group of storages. These mechanisms that do not require any individual (potentially private) information to exchanged among participants.  They are: i) adding an artificial term to the cost function of each storage unit; ii) longterm cooperation via an aggregating entity.  For each of these mechanisms, we show that self-interested storage units can obtain the profit obtained by a coalition of cooperating units. We show both mechanisms can achieve the maximum possible profit. 


This paper is organized as follows. In section \ref{section:Market and storage models} the market and storage models are introduced. We also introduce two scenarios: 1) the grand coalition solution where the aggregate profit is maximized and the Nash equilibrium where each storage plays a non-cooperative game. In section \ref{section: ACF} an approach to drive the profit of self-interested storages to the maximum possible profit via artificial cost functions is presented. In section \ref{section: Cooperation via aggregator} an aggregator model is introduced and the conditions for longterm cooperation between aggregator and storages are presented. In this same section, the Nash axiomatic bargaining model is presented and used to predict the profit split that the aggregator and storages would negotiate. Section \ref{section: Conclusion} concludes this paper.


\section{Market and storage models}
\label{section:Market and storage models}
In this section we describe the energy storage model and the market in which they interact. We also lay out two different scenarios: i) one where the storages cooperate to maximize the aggregate profit, and ii) another where the storages play a non-cooperative game and individually maximize profits.  We refer to the former as the grand coalition (GC) solution and to the latter as the Nash equilibrium (NE) solution. 

\subsection{Market model}
The price of energy at time $t$, $p^{[t]}\!\!\left(\boldsymbol{d}^{[t]}\right)$, is sensitive to energy demanded or supplied by the storages and is given by
\begin{align*}
p^{[t]}\!\!\left(\boldsymbol{d}^{[t]}\right)=\beta^{[t]} + \gamma^{[t]} \sum_{i \in \mathcal{I}}{ d_i^{[t]} } \quad \forall t \in \mathcal{T} 
\end{align*}
where $\boldsymbol{d}^{[t]}$ represents the actions of all storages at time $t$. Let $d_i^{[t]} \in \mathbb{R}$ denote the energy purchased (when $d_i^{[t]} > 0$) or sold (when $d_i^{[t]} < 0$) by the $i^\mathrm{th}$ storage at time $t$. The set of all storages is denoted by $\mathcal{I}$ and has size $n$. The constant $\beta^{[t]}$ is the price when the net purchases made by all storages is zero and $\gamma^{[t]}$ is a positive constant that determines the sensitivity of price to energy demand. The set of all time periods is denoted by $\mathcal{T}$ and has size $n_t$.

\subsection{Storage model}
Storages are agents that can buy energy at some time and sell it at another. The net energy purchased and sold of every storage is required to be zero and is expressed by
\begin{align}
\sum_{t \in \mathcal{T}}{d_i^{[t]}  } =0 \quad \forall i \in \mathcal{I}.\label{eq:net energy}
\end{align}
Because this paper seeks to emphasize the interaction between storages, other constraints such as energy and/or power limit are modeled by a cost function associated with each storage.

The profit for the $i^{\mathrm{th}}$ storage can be expressed as
\begin{align*}
\pi_i(\boldsymbol{d}_i;\boldsymbol{d}_{-i}) = \sum_{t \in \mathcal{T}}{  \left\{  -p^{[t]}\!\!\left(\boldsymbol{d}^{[t]}\right) \cdot d_i^{[t]}  { -  c_{i}\left(d_i^{[t]}\right)} \right\}} 
\end{align*}
where $\boldsymbol{d}_i $ and $\boldsymbol{d}_{-i}$  are the strategy choices of storage $i$ and the strategy choices of all storages excluding storage $i$, respectively. We assume a bounded strategy space for all storages. The battery degradation, efficiency, and/or energy transaction costs of storage $i$ are represented by the cost function $c_{i}(\cdot)$. 
It is known that as the depth of discharge increases, the costs of utilizing storage increases faster than linear \cite{Koller, Ortega-Vazquez}. Throughout this paper we assume a quadratic function of the form $c_{i}(x) = \frac{\epsilon_i}{2}x^2$ that captures faster-than-linear increasing costs. The positive constant $\epsilon_i$ is a storage specific cost coefficient. 

Now we define the GC and the NE solutions. In the GC solution
\begin{align}
& \boldsymbol{d}^*=\argmax_{\boldsymbol{d}} \sum_{i \in \mathcal{I}}{\pi_i(\boldsymbol{d}_i;\boldsymbol{d}_{-i})} \label{gc}\\ 
& \mathrm{s.t.\quad}  \sum_{t \in \mathcal{T}}{d_i^{[t]}  } =0 \quad (\lambda_i)\quad\forall i \in \mathcal{I} \nonumber
\end{align}
the \emph{aggregate profit} of the energy storages is maximized. The dual variables of the equality constraints are denoted by $\lambda_i$ and the GC solution is denoted by $\boldsymbol{d}^*$.

In the NE solution
\begin{align}
& \boldsymbol{d}_i'=\argmax_{\boldsymbol{d}_i} \pi_i(\boldsymbol{d}_i;\boldsymbol{d}_{-i}) \quad\; \forall i \in \mathcal{I}\label{eq:NSP}\\ 
& \mathrm{s.t.} \quad\sum_{t \in \mathcal{T}}{ d_i^{[t]}   } =0 \quad (\lambda_i) \nonumber
\end{align}
each storage maximizes its \emph{own profit} given the strategy choices of all other storages. The Nash equilibrium (NE) for storage $i$ is denoted by $\boldsymbol{d}_i'$. In the NE, no player has an incentive to deviate form his or her strategy.   

For readability and to convey intuition about the problem we consider a two-period case throughout the rest of this section and section \ref{section: ACF}. In subsection \ref{subsec: general} we generalize our results to $n_t$ time periods. 

\subsection{Solution to the two-period GC and NE solutions}

In this subsection it is shown that the NE yields a lower aggregate profit with respect to the GC strategy. This motivates the need to ``fix" the NE solution. For simplicity and without loss of generality, let $\beta^{[2]} - \beta^{[1]} =1$. From constraint \eqref{eq:net energy}, $d_i^{[1]}=-d_i^{[2]}$ for the two period case. A solution for the $i^\mathrm{th}$ storage is denoted by $d_i=d_i^{[1]}=-d_i^{[2]}$. 

\subsubsection*{Lemma 1}
As $n \to \infty$, the aggregate profit given by the GC solution increases and approaches a finite positive number while the aggregate profit under the NE solution approaches zero. 

\subsubsection{Numerical example}
Throughout the numerical examples $\gamma=1$, $n=2$, and $n_t=2$. Both storages have cost coefficients of $\epsilon_i = 1$. 

The GC solution is $d_i^*=\sfrac{1}{6}$ (\emph{i.e.} each storage charges $\sfrac{1}{6}$ in the first period and discharges the same amount in the second period). On the other hand, the NE is $d_i'=\sfrac{1}{5}$. The aggregate profit under the GC solution is $\sfrac{1}{6}$ while the aggregate profit under the NE is lower at $\sfrac{4}{25}$. 

\begin{figure}[h]
  \centering
   \includegraphics[width=0.5\textwidth]{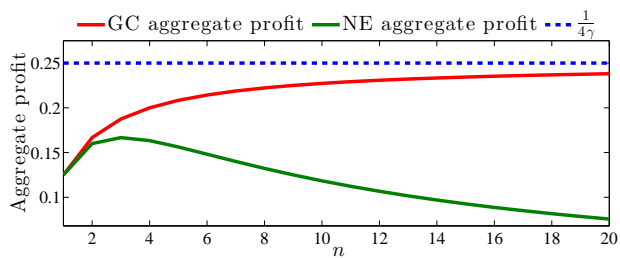}
      \caption{Aggregate profits for identical storages as a function of $n$ under the CG and NE solutions. The blue line is the limit of the GC profit as $n \rightarrow \infty$. The parameters are $\epsilon_i = 1\; \forall i \in \mathcal{I}$ and $\gamma=1$.} \label{Figure1}
\end{figure}

The NE oversupplies storage services with respect to the GC solution. The storages move more energy across time but the price difference between buying and selling periods is smaller and thus the NE profit is smaller. Figure \ref{Figure1} shows the aggregate profit for both the GC and the NE as a function of number of storages.

\subsubsection*{Proof of Lemma 1}
Using $d_i =d_i^{[1]}= -d_i^{[2]}$, $\beta^{[2]} - \beta^{[1]} = 1$, and $n_t=2$, the Lagrangian function of the problem \eqref{gc} is $\mathcal{L}(d_i) =  \sum_{i \in \mathcal{I}} \left\{d_i-d_i \gamma \sum_{j \in \mathcal{I}} d_i - \epsilon_id_i^{2}\right\}$ and its Karush-Kuhn-Tucker (KKT) optimality conditions
\begin{align}
&\frac{\partial\mathcal{L}(d_i) }{\partial d_i} =  1  - 2\gamma \sum_{j \in \mathcal{I}}d_j -2\epsilon_id_i = 0\quad   \forall\; i \in\; \mathcal{I}\label{eq:KKTGC_1}
\end{align}
are satisfied by $d_i^*= \sfrac{1}{2\epsilon_i \left( 1+ \gamma \sum_{j\in \mathcal{I}}{\frac{1}{\epsilon_j}}\right)}$. We show this by replacing $d_i$ in \eqref{eq:KKTGC_1} by $d_i^*$:
\begin{align*}
\frac{\partial\mathcal{L}(d_i^*) }{\partial d_i^*}\! \!&=  1-  \frac{\gamma\sum_{j \in \mathcal{I}}\frac{1}{\epsilon_j}}{ 1+\gamma \sum_{j\in\mathcal{I}} \frac{1}{\epsilon_j}} - \frac{1}{1+\gamma \sum_{j\in\mathcal{I}} \frac{1}{\epsilon_j}}=0\; \forall i\in \mathcal{I}.
\end{align*}

In the NE, no storage has the incentive to unilaterally change his or her strategy. Equivalently, for every $i$, $\boldsymbol{d}_i'$ solves $\max_{\sum_{t \in \mathcal{T}}{ d_i^{[t]}   } =0} \pi_i(\boldsymbol{d}_i;\boldsymbol{d}_{-i}')$ \cite{Gibbons}.
Each player's Lagrangian function of problem \eqref{eq:NSP} is $\mathcal{L}_i(d_i) = d_i -d_i  \gamma\sum_{j\in \mathcal{I}} d_j - \epsilon_id_i^2 \quad \forall i \in \mathcal{I} $ and their KKT optimality conditions
\begin{align}
\frac{\partial\mathcal{L}_i(d_i) }{\partial d_i} =  1 - \gamma \sum_{j \in \mathcal{I}}d_j -\left( \gamma+2\epsilon_i  \right) d_i = 0\quad  \forall\; i \in\; \mathcal{I} &\label{eq:KKTNE_1}
\end{align}
are satisfied by $d_i'=\sfrac{1}{\left(2\epsilon_i+\gamma\right) \left(1+\gamma \sum_{j \in \mathcal{I}} \frac{1}{2\epsilon_j+\gamma}\right)}\; \forall i \in \mathcal{I}$. When $\epsilon_i>0 \; \forall i \in \mathcal{I}$ or $\gamma > 0$, the system of equations described by \eqref{eq:KKTNE_1} has a unique solution. It follows that the NE is unique.  

We now show that the profit under the NE approaches zero while the profit under the GC approaches $\sfrac{1}{4\gamma}$ as the number of storages increases. The profit to be shared among the storages under the GC solution is
\begin{align*}
\sum_{i \in \mathcal{I}}{\!\! \pi_i(d_i^*,\boldsymbol{d}_{-i}^* ) }\!&=\!  \sum_{i \in \mathcal{I}}{\left\{d_i^*  \!-\! d_i^*\gamma \sum_{j \in \mathcal{I}}d_j^* \!-\! \epsilon_i d_i^{*2} \right\}} \\
& = \frac{\gamma\left(\sum_{j \in \mathcal{I}}\frac{1}{\epsilon_j} \right)^2}{4\left( 1+ \gamma\sum_{j \in \mathcal{I}}\frac{1}{\epsilon_j} \right)^2}. 
\end{align*}
As $n \to \infty$, $\sum_{j \in \mathcal{I}}\frac{1}{\epsilon_j} \to \infty$ and the aggregate profit made by the storages under the GC solution approaches $\sfrac{1}{4\gamma}$.

To show that the NE aggregate profit goes to zero as the number of storages increase, we assume $\epsilon_i=0\;\forall \;i \in \mathcal{I}$. Then, an upper bound for aggregate profit under NE is given by
\begin{align*}
\sum_{i \in \mathcal{I}}{\!\! \pi_i(d_i',\boldsymbol{d}_{-i}' ) }\!&=\!  \sum_{i \in \mathcal{I}}{\left\{d_i'  \!-\! d_i'\gamma \sum_{j \in \mathcal{I}}d_j'\right\}} \\
& = \frac{n}{\left( 1+ n\right)^2} \to 0\; \mathrm{as}\; n \to \infty.
\end{align*}
The lower bound must be non-negative as he or she can always choose $d_i=0$ to achieve a zero profit. It follows that as $n \to \infty$, the aggregate profit under the NE approaches zero.  \QEDB


%
%

\section{Fixing the Nash equilibrium via Artificial Cost Functions}
\label{section: ACF}
In this section we study the use of artificial cost functions (ACFs) to ``fix" the NE. We would like to find a set of ACFs $g_i^{[t]}(\cdot)$ such that when the storages incur it, the NE is equal to the GC solution. The ACF is effectively a control signal that penalizes deviations from the GC solution. We refer to the NE under the ACF as the ``artificial" Nash equilibrium (ANE). Having a NE that equals the GC solution is desirable because: 1) the GC aggregate profit is larger than the NE aggregate profit, 2) it is strategically stable, and 3) it is self-enforcing \cite{Gibbons}. 

The idea of fixing an undesirable NE outcomes using a cost function was presented in \cite{Maheswaran}. However, revenue neutrality (which we will define shortly) is not a concern in their context. 
\subsection{Nash equilibrium under artificial cost functions}
The two-period profit for storage $i$ when exposed to the ACF is
\begin{align*}
\pi_i^\mathrm{A}(\boldsymbol{d}_i;\boldsymbol{d}_{-i}) &= d_i -d_i \gamma \sum_{i \in \mathcal{I}} d_i -\frac{\epsilon_i}{2}d_i^2 - g_i(d_i)
\end{align*}
where $g_i(d_i)=g_i^{[1]}(d_i)+g_i^{[2]}(-d_i)$ and the ANE solution is denoted by 
\begin{align}
& \overline{d}_i =\argmax_{d_i} \pi_i^\mathrm{A}(\boldsymbol{d}_i;\boldsymbol{d}_{-i}) \quad\; \forall i \in \mathcal{I}.\label{artificial_NE}
\end{align}

\subsubsection*{Lemma 2}
There exists a cost function of the form $g_i(d_i)=a_id_i^{2}+ b_id_i \; \forall i \in \mathcal{I},\; t \in \mathcal{T}$ such that:
\begin{itemize}
\item The ANE solution of problems \eqref{artificial_NE} equals the solution of problem \eqref{gc} (\emph{i.e.} $\overline{d}_i=d_i^*\; \forall i \in \mathcal{I}$).
\item  It is revenue neutral (\emph{i.e.} $g_i(\overline{d}_i)=0 \;\forall i \in \mathcal{I}$).
\end{itemize}

Moreover, the coefficients of $g_i(\cdot)$ are given by
\begin{align*}
a_i \!&= -\gamma\left( 1- \epsilon_i \sum_{j \in \mathcal{I}}{\frac{1}{\epsilon_j}}\right) \quad\forall i \in \mathcal{I} \\
 b_i\!&= -\frac{a_i}{2\epsilon_i\left( 1+ \gamma \sum_{j \in \mathcal{I}}{\frac{1}{\epsilon_j}}\right)} \quad \forall i \in \mathcal{I}.
 \end{align*}
Note that $a_i$ and $b_i$ only depend on individual information ($\epsilon_i$), public information ($\gamma$), and on the sum of other storages' characteristics ($\sum_{j \in \mathcal{I}}{ \! \frac{1}{\epsilon_j } }$). The implication of this is that the GC solution can be reached in a distributed fashion, without the need of each storage disclosing its information to the rest of the storages. 

\subsubsection*{Proof of Lemma 2}
The coefficients of the ACF and $d_i$ must satisfy the following system of equations
\begin{subequations}
\begin{align}
\frac{\partial\pi_i^A(\boldsymbol{d}_i;\!\boldsymbol{d}_{-i})}{\partial d_i} \!\!&= \!\!1\!-\!b_i \!-\! \gamma d_i \!-\! \gamma\! \!\sum_{j\in \mathcal{I}}{\!\!d_j} \!- \!2(\!\epsilon_i \!+\!a_i\!) d_i \!= \!0 \; \forall i \in \mathcal{I}. \label{ACF system 1}\\
\frac{\partial\mathcal{L}(d_i) }{\partial d_i} &=  1  - 2\gamma \sum_{j \in \mathcal{I}}d_j -2\epsilon_id_i = 0\quad   \forall\; i \in\; \mathcal{I} \label{ACF system 2} \\
 a_i d_i &= - b_i \qquad \forall i \in \mathcal{I}.\label{ACF system 3}
\end{align}
\end{subequations}
Equations \eqref{ACF system 1} and \eqref{ACF system 2} ensure that in addition to satisfying each player's individual profit maximization problem, the ANE satisfies the GC solution. Equation \eqref{ACF system 3} ensures revenue neutrality. From $\beta^{[2]}- \beta^{[1]}-=1$, the solutions of both the GC and the NE are non-negative. Thus $a_i d_i^2 +  b_id_i=0$ is replaced by \eqref{ACF system 3}. 

It is straight forward to show that $a_i$, $b_i$, and $d_i^*$ satisfy \eqref{ACF system 1} by using the expressions shown in Lemmas 1 and 2
\begin{align*}
\frac{\partial\pi_i^A(\boldsymbol{d}_i;\!\boldsymbol{d}_{-i})}{\partial d_i} \!\!&= \!\!1\! - \! \frac{\gamma\!\left( \!1\!-\! \epsilon_i\! \sum_{j \in \mathcal{I}}{\frac{1}{\epsilon_j}}\right) \!-\! \gamma \epsilon_i \!\sum_{j \in \mathcal{I}}{\frac{1}{\epsilon_j}}\!-\!  z_i }{2\epsilon_i\left( 1+ \gamma \sum_{j \in \mathcal{I}}{\frac{1}{\epsilon_j}}\right)} =0 
\end{align*}
where $z_i = 2\left(\epsilon + a_i \right) + \gamma$. By Lemma 1, $d_i^*$ satisfies \eqref{ACF system 2}. We can conclude that $\overline{d}_i=d_i^*$. 

Finally, we show revenue neutrality (\emph{i.e.} $g_i(\overline{d}_i) = 0 \; \forall i \in \mathcal{I}$): 
\begin{align*}
g_i(\overline{d}_i) \!\!&= \!a_i\overline{d}_i^2 + b_i\overline{d}_i \\
&= \!\frac{-\gamma \left(1-\epsilon_i \sum_{j \in \mathcal{I} } \frac{1}{\epsilon_j} \right) }{4\epsilon_i^2\! \left(1+\gamma\sum_{j \in \mathcal{I}}{\frac{1}{\epsilon_j} }\right)^2} \!+\! \frac{\gamma \left(1-\epsilon_i \sum_{j \in \mathcal{I} } \frac{1}{\epsilon_j} \right) }{4\epsilon_i^2 \!\left(1+\gamma\sum_{j \in \mathcal{I}}{\frac{1}{\epsilon_j} }\right)^2}\!=\!0.
\end{align*}

\QEDB

\subsection{Sensitivity analysis of the artificial cost function} 
\begin{figure}[h!]
  \centering
    \includegraphics[width=0.5\textwidth]{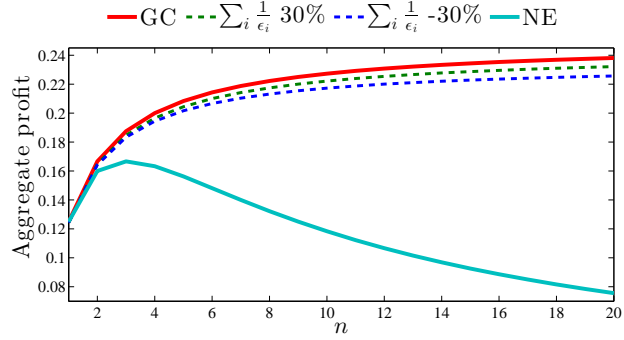}
      \caption{ Aggregate profits as a function of number of storages participating. The lines $\sum_{i} \frac{1}{\epsilon_i}\; 30\%$ and $\sum_{i} \frac{1}{\epsilon_i}\; -30\%$ correspond to overestimation and underestimation, respectively, of $\sum_{i \in \mathcal{I}} \frac{1}{\epsilon_i}$.   }
      \label{sensitivity}
\end{figure}

\begin{figure}[h!]
  \centering
    \includegraphics[width=0.5\textwidth]{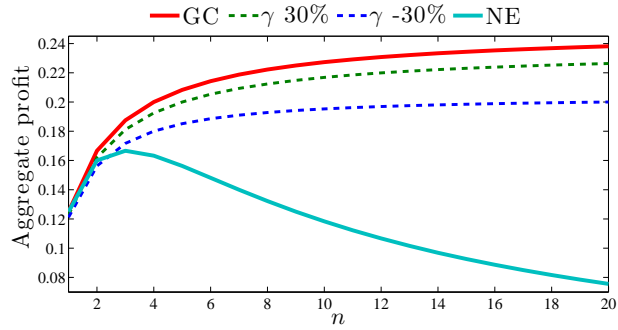}
      \caption{ Aggregate profits as a function of number of storages participating. The lines $\gamma\; 30\%$ and $\gamma\; -30\%$ correspond to overestimation and underestimation, respectively, of $\gamma$.   }
      \label{sensitivity2}
\end{figure}

In this subsection we show the ACF aggregate profit is robust to misestimations of the parameters needed to compute the ACF.  Figures \ref{sensitivity} and \ref{sensitivity2} show the effect of overestimating or underestimating $\sum_{i \in \mathcal{I}} \frac{1}{\epsilon_i}$  and $\gamma$, respectively, by 30\% on the aggregate profit. Even with large misestimations, the GC aggregate profit remains considerably higher than the NE aggregate profit for most $n$. 

\subsection{Generalization to $n_t$ time periods}
\label{subsec: general}
In this subsection we generalize Lemma 2 to an arbitrary number of periods. 

\subsubsection*{Lemma 3}
There exist a set cost functions of the form $g_i^{[t]}(d_i^{[t]})=\frac{a_i^{[t]}}{2}d_i^{[t]2}+ b_i^{[t]}d_i^{[t]} \; \forall i \in \mathcal{I},\; t \in \mathcal{T}$ such that:
\begin{itemize}
\item The ANE of problems \eqref{eq:NSP} equals the solution of problem \eqref{gc} (\emph{i.e.} $\overline{d}_i^{[t]}=d_i^{[t]*}\; \forall i \in \mathcal{I},\; t \in \mathcal{T}$).
\item  It is revenue neutral (\emph{i.e.} $\sum_{t \in \mathcal{T}}g_i^{[t]}(\overline{d}_i^{[t]})=0 \;\forall i \in \mathcal{I}$).
\end{itemize}

\begin{figure}[h!]
  \centering
    \includegraphics[width=0.5\textwidth]{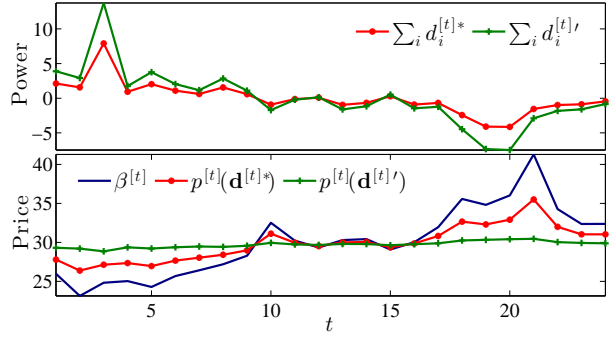}
      \caption{  The upper plot in Figure \ref{n_t_time_periods} shows the total energy purchases under both the NE $\sum_{i \in \mathcal{I}} {d_i^{[t]}}'$ and GC $\sum_{i \in \mathcal{I}} d_i^{[t]*}$. The lower plot shows the parameter $\beta^{[t]}$, price under the GC solution $p^{[t]}\left(\boldsymbol{d}^{[t]*}\right)$, and price under the NE solution $p^{[t]}\left({\boldsymbol{d}^{[t]}}'\right)$. The price parameter $\beta^{[t]}$ is the day-ahead price in the PJM interconnection during 01/02/2011 \cite{PJM}. The parameter $\gamma^{[t]}$ is randomly generated but proportional to $\beta^{[t]}$ to simulate the increasing slope of the typical energy supply curve. The number of storages is $n=20$ and their cost coefficients are randomly generated. }
      \label{n_t_time_periods}
\end{figure}

The upper plot in Figure \ref{n_t_time_periods} shows the total energy purchases under both the NE and GC. The lower plot in Figure \ref{n_t_time_periods} shows the price under both the NE and GC in a 24 time period game. Because more energy is moved in the NE solution, the price under the NE solution, $p^{[t]}\left({\boldsymbol{d}^{[t]}}'\right)$, is considerably flatter than under the GC solution $p^{[t]}\left(\boldsymbol{d}^{[t]*}\right)$. Even though the NE solution moves more energy across time, the price difference is smaller and their profits are lower. For this reason the storages may be inclined to cooperate to increase their aggregate profit. 

\subsubsection*{Proof of Lemma 3}
We would like to show the existence of the set of functions $g_i^{[t]}(\cdot)$. To do so, we restrict ourselves to the subset of functions whose revenue is zero for every period (\emph{i.e.} $g_i^{[t]}(\overline{d}_i^{[t]})=0 \;\forall i \in \mathcal{I},\;t \in \mathcal{T}$). Similar to the proof of Lemma 2, in order to find the coefficients of the artificial cost functions, we solve the following system of equations for  $d_i^{[t]}$, $a_i^{[t]}$,  $b_i^{[t]}\; \forall i \in \mathcal{I},\; t \in \mathcal{T}$ and $\lambda_i\; \forall i \in \mathcal{I}$
\begin{subequations}
\label{eq: system 1}
\begin{align}
&\frac{\partial\mathcal{L}_i^A(d_i^{[t]}\! , \lambda_i) }{\partial d_i^{[t]}} \!\!=\!\!  - \beta^{[t]}\! - b_i^{[t]}\!   -\! \gamma^{[t]} \!\sum_{j \in \mathcal{I}}\!d_j^{[t]} \!-\!\left(\!\epsilon_i  \!+\! \gamma^{[t]}\!+a_i^{[t]}\!\right)\!d_i^{[t]} \nonumber\\
& \hspace{100 pt}+\lambda_i \!= \!0\quad \forall\; i \in\; \mathcal{I},\; t\in\;\mathcal{T} \label{eq:KKTACF_T_1}\\
&\frac{\partial\mathcal{L}(d_i^{[t]}, \lambda_i) }{\partial d_i^{[t]}} =  - \beta^{[t]}  - 2\gamma^{[t]} \sum_{j \in \mathcal{I}}d_j^{[t]} -\epsilon_id_i^{[t]} +\lambda_i = 0 \label{eq:KKTGC_T_1}\\
& \hspace{145 pt}  \forall\; i \in\; \mathcal{I},\; t\in\;\mathcal{T} \nonumber\\
&\frac{\partial\mathcal{L}(d_i^{[t]}, \lambda_i) }{\partial \lambda_i^{[t]}} =  \sum_{t \in \mathcal{T}}{d_i^{[t]}}  =  0 \quad \forall\; i \in\; \mathcal{I} &\label{eq:EQC_T_1} \\
& \frac{a_i^{[t]}}{2}d_i^{[t]2}+ b_i^{[t]}d_i^{[t]} \; \forall i \in \mathcal{I},\; t \in \mathcal{T}\label{eq:RN_T_1} 
\end{align}
\end{subequations}
where $\mathcal{L}_i^A(d_i^{[t]}\! , \lambda_i)$ is the Lagrangian function of the $n_t$ time periods individual profit maximization problem. Equations \eqref{eq:KKTACF_T_1} ensure that $d_i^{[t]}$ satisfies the ANE. Equations \eqref{eq:KKTGC_T_1} ensure that the solution satisfies the GC solution while equations \eqref{eq:EQC_T_1} enforce the equality constraints of each storage. Finally, equation \eqref{eq:RN_T_1} ensures revenue neutrality. 

The solution to the multiple period GC problem is $d_i^{[t]*}=\frac{\sum_k \frac{\beta^{[k]}-\beta^{[t]}}{z^{[k]}} }{ \epsilon_iz^{[t]}\sum_k\frac{1}{z^{[k]}}} \; \forall\; i \in\; \mathcal{I},\; t\in\;\mathcal{T}$ and $\lambda_i^* =  \frac{\sum_t \frac{\beta^{[t]}}{z^{[t]}}}{\sum_t\frac{1}{z^{[t]}}} \; \forall\; i \in\; \mathcal{I}$ where $z^{[t]}=1+2\gamma^{[t]} \sum_j\epsilon_j^{-1}$. It is straight forward to show that they satisfy the KKT conditions of the GC problem given by equations \eqref{eq:KKTACF_T_1} and  \eqref{eq:EQC_T_1}. 

From the multiple period GC solution, when $\beta^{[t]} < \frac{\sum_k \frac{\beta^{[k]}}{z^{[k]}} }{ \sum_k\frac{1}{z^{[k]}}}$ then $d_i^{[t]*}> 0\; \forall\; i \in\; \mathcal{I},\; t\in\;\mathcal{T}$. Denote the set of such time periods as $\mathcal{T}_1$. The rest of the time periods (when $d_i^{[t]}$ is non-positive) are in the set $\mathcal{T}_2$. We can then replace the ACF revenue neutrality requirement \eqref{eq:RN_T_1} with
\begin{subequations}
\label{eq: RN_T_2}
\begin{align}
&a_i^{[t]} d_i^{[t]} +2 b_i^{[t]} =0 \quad\forall \;i \in \mathcal{I},\; \forall t \in \mathcal{T}_1 \label{eq:revneut_1}\\  
&a_i^{[t]} d_i^{[t]} - 2 b_i^{[t]} = 0 \quad\forall \;i \in \mathcal{I},\; \forall t \in \mathcal{T}_2.\label{eq:revneut_2}
\end{align}
\end{subequations}

We substitute the term $a_i^{[t]}d_i^{[t]}$ in \eqref{eq:KKTACF_T_1} with either $-2 b_i^{[t]}$ or $2 b_i^{[t]}$ depending on whether $t$ is in $\mathcal{T}_1$ or in $\mathcal{T}_2$. Finally, instead of solving the system of non-linear equations described by \eqref{eq: system 1} we solve the following system of linear equations for the variables $d_i^{[t]}$, $\lambda_i$, and $b_i^{[t]}$:
\begin{subequations}
\begin{align}
&  - \beta^{[t]}\! + b_i^{[t]}\!   -\! \gamma^{[t]} \!\sum_{j \in \mathcal{I}}\!d_j^{[t]} \!\!-\!\!\left(\!\epsilon_i  \!+\! \gamma^{[t]}\!\right)\!d_i^{[t]} \! \!+\!\!\lambda_i \!= \!0\;\forall\; i \!\in\! \mathcal{I}, t\!\in\!\mathcal{T}_1 \label{eq: system_2_1} \\
&- \beta^{[t]}\! -3 b_i^{[t]}\!   -\! \gamma^{[t]} \!\sum_{j \in \mathcal{I}}\!d_j^{[t]} \!\!-\!\!\left(\!\epsilon_i  \!+\! \gamma^{[t]}\!\right)\!d_i^{[t]} \! \!+\!\!\lambda_i \!= \!0\;\forall\; i \!\in\! \mathcal{I}, t\!\in\!\mathcal{T}_2 \label{eq: system_2_2}  \\
& \eqref{eq:KKTGC_T_1},\; \eqref{eq:EQC_T_1}. \nonumber
\end{align}
\end{subequations}

Equations \eqref{eq: system_2_1}, \eqref{eq: system_2_2}, and \eqref{eq:KKTGC_T_1} can be expressed in matrix notation as
\begin{align*}
& \matr{M}_j^{[t]} \begin{bmatrix} \boldsymbol{d}^{[t]} \\ 
 \boldsymbol{b}^{[t]} \end{bmatrix} 
 +  \matr{I}_{2} \boldsymbol{\lambda} = 
  \begin{bmatrix} \boldsymbol{1} \beta^{[t]}
  \\ \boldsymbol{1} \beta^{[t]}
  \end{bmatrix} \quad \forall\; t \in \mathcal{T}_j,\; j=1,\;2
  \end{align*}
  where
  \begin{align*}
  & \matr{M}_j^{[t]}=\begin{bmatrix} 
-\gamma^{[t]}\boldsymbol{1}\boldsymbol{1}^T - \matr{E} -\gamma^{[t]}\matr{I}_1 & c_j \matr{I}_1 \\
-2\gamma^{[t]}\boldsymbol{1}\boldsymbol{1}^T - \matr{E}& \matr{0}
\end {bmatrix},\; c_1 = 1,\; c_2=-3.
\end{align*}
 The vector $\boldsymbol{1}\in \mathbb{R}^n$ is an all ones vector and $\matr{E} \in \mathbb{R}^{n\times n} $ is a diagonal matrix whose $ii^\mathrm{th}$ entry is $\epsilon_i$. The symbol $\matr{I}_j \in \mathbb{R}^{j \cdot n\times n}$ represents $j$ vertically concatenated identity matrices. The $i^\mathrm{th}$ entries of vectors $\boldsymbol{d}^{[t]}\in \mathbb{R}^n$, $\boldsymbol{b}^{[t]}\in \mathbb{R}^n$, and $\boldsymbol{\lambda}\in \mathbb{R}^n$ are $d_i^{[t]}$, $b_i^{[t]}$, and $\lambda_i$ respectively.
 
 We can further compact equations \eqref{eq: system_2_1}, \eqref{eq: system_2_2}, \eqref{eq:KKTGC_T_1}, and \eqref{eq:EQC_T_1} to
 \begin{align*}
 \begin{bmatrix}\matr{M} & \matr{I}_{2 n_t} \\ \matr{N} &  \matr{Z} \end{bmatrix} \begin{bmatrix} \boldsymbol{d} \\ 
 \boldsymbol{b} \\ \boldsymbol{\lambda} \end{bmatrix} = \begin{bmatrix} \boldsymbol{\beta}
  \\ \boldsymbol{ \beta} \\ \boldsymbol{0} \end{bmatrix}
 \end{align*}
 where $\matr{M} \in \mathbb{R}^{2n\cdot n_t \times 2n\cdot n_t}$ is a block diagonal matrix whose $tt^\mathrm{th}$ block is $\matr{M}_i^{[t]}$, $\boldsymbol{d}=[\boldsymbol{d}^{[1]T},\hdots\boldsymbol{d}^{[n_t]T}]^T$, $\matr{M}_j^{[t]}$, $\boldsymbol{b}=[\boldsymbol{b}^{[1]T}\hdots\boldsymbol{b}^{[n_t]T}]^T$, and $\boldsymbol{\beta}=[\beta^{[1]}\boldsymbol{1}^T\hdots\beta^{[n_t]}\boldsymbol{1}^T]^T$. 
 The equation $\begin{bmatrix}\matr{N} & \matr{Z} \end{bmatrix}\begin{bmatrix} \boldsymbol{d}^T &
 \boldsymbol{b}^T & \boldsymbol{\lambda}^T \end{bmatrix}^T = \boldsymbol{0}$ represents equations \eqref{eq:EQC_T_1}. The matrix $\matr{N} \in \mathbb{R}^{ n \times 2 n_t \cdot n}$ is constructed by horizontally concatenating $\begin{bmatrix}\matr{I}_1 & \matr{Z} \end{bmatrix}$ $n_t$ times where $\matr{Z} \in \mathbb{R}^{n \times n}$ all-zero matrix.
 
 Using Gaussian elimination, it is straightforward to show that matrix $ \begin{bmatrix}\matr{M} & \matr{I}_{2\cdot n_t} \\ \matr{N} &  \matr{Z} \end{bmatrix} $ is full rank. It follows that it is invertible and the system of equations \eqref{eq: system_2_1}, \eqref{eq: system_2_2}, \eqref{eq:KKTGC_T_1}, and \eqref{eq:EQC_T_1} has a unique solution.
 
Equations \eqref{eq:revneut_1} and \eqref{eq:revneut_2} can be used to find coefficients $a_i^{[t]}$. \QEDB


\section{Cooperation via Aggregator} 
\label{section: Cooperation via aggregator}
In this section we explore the possibility of the storages reaching the GC solution by cooperating with a central entity that we refer to as the ``aggregator.'' In this setting, the storages do not have access to the wholesale market but instead buy/sell energy from/to an aggregator. The aggregator determines the prices paid by/to the storages. 

Previous work on aggregators \cite{Ortega-Vazquez_TPS_May_2005, Wu, DiWu} assume cooperation between aggregator and storages. In this work, however, we analyze possible outcomes of the aggregator-storage interaction and do not assume that the aggregator will cooperate with the storages and vice versa. The aggregator-storage game is modeled as a simultaneous move game. We show that the single-shot NE is inefficient, explore the possibility of aggregator-storage cooperation in the longterm, and derive conditions for cooperation.  

\subsection{Aggregator model}
The aggregator profits by purchasing or selling energy on the wholesale market described in section \ref{section:Market and storage models} and in turn selling to or buying from the storages. For every time period, the aggregator determines the price of energy that each storage pays and the storages decide how much to purchase from or sell to the aggregator. The prices sent by the aggregator are assumed to be bounded. The aggregator's profit from trading with player $i$ is denoted by
\begin{align*}
\pi_{a,i}\!\left(\boldsymbol{\tau}_i; \boldsymbol{d}_i\right) \!=  \sum_{t\in\mathcal{T}} \left\{\tau^{[t]}_i d_i^{[t]} -  p^{[t]}\!\!\left(\boldsymbol{d}^{[t]}\right) \cdot d_i^{[t]}  \right\} 
\end{align*}
where $\boldsymbol{\tau}_i$ and $\boldsymbol{d}_i$ are the strategies of the aggregator and storage $i$, respectively. The strategy space of the aggregator is the set of possible price schedules that it can send to the storage. The strategy space of each storage is the set of all feasible charge/discharge schedules and is assumed to be bounded. The energy price that storage $i$ pays at time $t$ is denoted by $\tau_i^{[t]}$.

\subsection{Storage problem under an aggregator}
With the aggregator acting as a middle-man between the wholesale market and the storages, the storages are insensitive to the wholesale market prices, and only respond to the prices sent by the aggregator. Let
 \begin{align*}
\pi_i(\boldsymbol{d}_i; \boldsymbol{\tau}_i) =\sum_{t \in \mathcal{T}}\left\{- \tau_i^{[t]}d_i^{[t]} - \frac{\epsilon_i}{2}d_i^{[t]2} \right\}
\end{align*}
denote the profit of the $i^\mathrm{th}$ storage under an aggregator.

In the NE solution under an aggregator, the storage and aggregator make their strategy choices either simultaneously or without knowledge of the other player's choices. The game can be expressed as 
\begin{subequations}
\label{so_storage_game}
\begin{align}
& \doverline{\boldsymbol{\tau}}_i=\argmax_{\tau_i^{\text{min}} \le \tau_i^{[t]} \le \tau_i^{\text{max}} } \pi_{a,i}\!\left(\boldsymbol{\tau}_i; \boldsymbol{d}_i\right) \label{so_problem} \\ 
&\doverline{\boldsymbol{d}}_i=\argmax_{\sum_{t\in\mathcal{T}}{d_i^{[t]}}=0 \quad \left(\lambda_i\right)}\pi_i(\boldsymbol{d}_i; \boldsymbol{\tau}_i) 
\label{eq:storage_problem}
\end{align}
\end{subequations}
 where every player independently maximizes its own profit. The prices that the aggregator can send to the storage are bounded by $\tau_i^{\text{min}}$ and $\tau_i^{\text{max}}$.

\subsubsection*{Lemma 4}
Both players earn a profit of zero in all Nash equilibria of game \eqref{so_storage_game}.
\subsubsection{Numerical example}
Assume that $\tau_i^{\text{min}} = -1$ and $ \tau_i^{\text{max}} =1$. Denote $\tau_i = \tau_i^{[2]} - \tau_i^{[1]}$. 
First we analyze the aggregator's best response given the storage's strategy. If the storage chooses $d_i>0$, the aggregator's best response is $\tau_i = -2$. On the other hand if the storage chooses $d_i<0$,the aggregators best response is $\tau_i = 2$. If the storage chooses $d_i=0$, the aggregator's best response is any feasible $\tau_i$. 

Now we analyze the storage's best response given the aggregator's strategy. Both storages maximize their profits by playing $d_i =\sfrac{\tau_i}{2}$. It is easy to see that the only stable NE is $\doverline{d}_i=0$ and $\doverline{\tau}_i$ is any feasible $\tau_i$. It follows that the profit is zero for both the aggregator and the storage as no transactions occur. 

 It is likely that the aggregator and storages interact repeatedly over time. Thus, there might be the possibility of fostering longterm cooperation between the storages and the aggregator to achieve the GC solution. We will present a infinitely repeated game model and show under which conditions cooperation can be sustained in the longterm. 
\subsubsection*{Proof of Lemma 4}
It is straight forward to show that, the storages actions, the aggregator maximizes its profit by sending price schedule
\begin{align*}
\doverline{\tau}_i^{[t]}=\begin{cases} 
\tau_i^{\text{min}} \;\; \mathrm{if}\;\; d_i^{[t]}<0  \\
\tau_i^{\text{max}} \;\; \mathrm{if}\;\; d_i^{[t]}>0  \\
\tau_i^{\text{min}} \le \tau_i^{[t]}\le \tau_i^{\text{max}} \;\; \mathrm{if}\;\; d_i^{[t]}=0  
\end{cases} \forall  \;t \in \mathcal{T}.
\end{align*} 

The KKT conditions of problem \eqref{eq:storage_problem} are 
\begin{align*}
-\tau^{[t]} - \epsilon_id_i^{[t]} +\lambda_i&= 0 \quad \forall t \in \mathcal{T} \\
 \sum_{t \in \mathcal{T}} d_i^{[t]} &= 0  
\end{align*}
and are satisfied by $\doverline{d}_i^{[t]} = \frac{\frac{1}{n_t} \sum_{k\in \mathcal{T}}\tau_i^{[k]}  - \tau_i^{[t]}}{\epsilon_i} \; \forall t\in\mathcal{T}$ for a given set of $\tau_i^{[t]}$ and $\doverline{\lambda}_i=\frac{1}{n_t} \sum_{k \in \mathcal{T}}\tau_i^{[k]} $.

The Nash equilibria should satisfy problems \eqref{so_storage_game}. It follows the Nash equilibria is $\tau_i^{[t]}=a \;\forall \;t \in \mathcal{T}$, where $a$ is a constant such that $\tau_i^{\mathrm{min}}\le a \le \tau_i^\mathrm{max}$, and $d_i^{[t]}=0 \;\forall\;t \in \mathcal{T}$. Any other strategy choices are unstable. Since the storage does not purchase or sell energy in all Nash equilibria, the profit for both the aggregator and storage is zero. \QEDB
\subsection{Repeated game model}
 We assume that the single-stage game \eqref{so_storage_game} is repeated indefinitely. The longterm profit made by each player is the discounted sum of the single-stage profits. Denote the longterm profit made by the aggregator from trading with storage $i$ as 
\begin{align*}
\pi^{\infty}_{a,i} = \sum_{k=0}^{\infty}{ \delta^k{\pi_{a,i}(\boldsymbol{\tau}_i^{(\!k\!)}; \boldsymbol{d}_i^{(\!k\!)})}  } 
\end{align*}
where $ \delta \in (0,1)$ is the discount rate (\emph{e.g.} interest rate). The symbols $\boldsymbol{\tau}_i^{(\!k\!)}$ and $\boldsymbol{d}_i^{(\!k\!)}$ denote strategy decisions for the $k^\mathrm{th}$ time the single-stage game is played. 
Similarly, the longterm profit of storage $i$ is given by
 \begin{align*} 
 \pi^{\infty}_i = \sum_{k=0}^{\infty}{ \delta^k\pi_i(\boldsymbol{d}_i^{(\!k\!)}; \boldsymbol{\tau}_i^{(\!k\!)})  } . 
 \end{align*}

\subsubsection{Strategy space for the repeated game}
\label{coop_strategies}
In order to keep the repeated game tractable, the strategy spaces of the aggregator and storages are reduced to specific cooperation and defection strategies.
\subsubsection{Cooperation strategies}
The cooperation strategy of the aggregator is given by
\begin{align*}
\boldsymbol{\tau}_i^{(\!k\!)}=\begin{cases} 
\hat{\boldsymbol{\tau}}_i \;\; \mathrm{if}\;\; \boldsymbol{d}_i^{(m)}=\hat{\boldsymbol{d}}_i \quad \forall \; m <k \\
\doverline{\boldsymbol{\tau}}_i \;\; \text{otherwise}
\end{cases} \forall i \in \mathcal{I}.
\end{align*} 
It describes the strategy in which during the $k^\mathrm{th}$ game, the aggregator sends storage $i$ a previously agreed price schedule $\hat{\boldsymbol{\tau}}_i$ if storage $i$ has played an agreed $\hat{\boldsymbol{d}}_i$ during all previous times. If storage $i$ fails to uphold its commitment, the aggregator stops cooperating and plays the NE solution ($\doverline{\boldsymbol{\tau}}_i$) for the subsequent times the game is played. Likewise, 
\begin{align*}
\boldsymbol{d}_i^{(\!k\!)}=\begin{cases} 
\hat{\boldsymbol{d}}_i \;\; \mathrm{if}\;\; \boldsymbol{\tau}_i^{(m)}=\hat{\boldsymbol{\tau}}_i \quad \forall \; m<k \\
\doverline{\boldsymbol{d}}_i \;\; \text{otherwise}
\end{cases} \forall i \in \mathcal{I} 
\end{align*}
describes the cooperation strategy of storage $i$. Storage $i$ plays an agreed $\hat{\boldsymbol{d}}_i$ if the aggregator has upheld its commitment to send an agreed $\hat{\boldsymbol{\tau}}_i$ during all previous times the game has been played. If the aggregator fails to uphold its commitment, the storage $i$ stops cooperating and plays the NE solution ($\doverline{\boldsymbol{d}}_i$) for the subsequent times the game is played.

\subsubsection{Defection strategies} 
We now describe the defection or ``cheating" strategies that the aggregator and the storage can play. We assume that all players are aware that, as described in section \ref{coop_strategies}, players stop cooperating when the opponent fails to uphold its commitment. Therefore, if one of the players decides to cheat, it will do so by maximizing its single-stage profit. The aggregator's single-stage profit derived from storage $i$ is maximized during a single game by playing 
\begin{align*}
\boldsymbol{\tau}_i^D = \argmax_{\tau^{\text{min}}_i \le \tau_i^{[t]} \le \tau^{\text{max}}_i  } \pi_{a,i}(\boldsymbol{\tau}_i; \hat{\boldsymbol{d}}_i).  
\end{align*}
Similarly, storage $i$ maximizes its profit for a given $\hat{\boldsymbol{\tau}}_i$ during a single game by playing the defection strategy 
\begin{align*}
\boldsymbol{d}_i^D = \argmax_{\sum_id_i^{[t]} = 0}\pi_i(\boldsymbol{d}_i; \hat{\boldsymbol{\tau}}_i)  \quad \forall i \in \mathcal{I}. 
\end{align*}

In the following subsection we show, given $\hat{\boldsymbol{d}}_i$, which choices of $\hat{\boldsymbol{\tau}}_i$ ensure that every player never plays its defection strategy.
\subsection{Ensuring cooperation in an infinitely repeated game} 
We would like to choose a $\hat{\boldsymbol{\tau}}_i$ (or equivalently, a profit split between storages and aggregator) such that cooperation is sustained by all players. Problems 
\begin{align*}
&v_i^*= \argsup_{v_i \in \mathbb{R}_+} \pi^{\infty}_{a,i} \;\; \mathrm{and} \;\; w_i^* = \argsup_{w_i \in \mathbb{R}_+} \pi^{\infty}_i \quad \forall i \in \mathcal{I} \nonumber
\end{align*}
are solved by the aggregator and the storages, respectively, to determine when to defect (if at all). The strategy $v_i \in \mathbb{R}_+$ is the time the aggregator decides to defect from cooperation with storage $i$. Likewise, $w_i \in \mathbb{R}_+$ denotes the time storage $i$ decides to stop cooperating with the aggregator. To ensure longterm cooperation, $w_i^* = \infty \;\forall i \in \mathcal{I}$ and $ v_i^* = \infty \; \forall i \in \mathcal{I}$. 

\subsubsection*{Lemma 5}

Cooperation with every storage is sustained by the aggregator, or equivalently, $v_i^*=\infty \; \forall i \in \mathcal{I}$ when all $\hat{\boldsymbol{\tau}}_i$ are in the sets $\mathcal{A}_i^{\mathrm{a}}  = \{ \boldsymbol{\tau}_i | \pi_{a}(\boldsymbol{\tau}_{n:n}; \hat{\boldsymbol{d}}_{n:n}  )\ge(1\!-\!\delta)\pi_{a}(\boldsymbol{\tau}^{D_n}; \hat{\boldsymbol{d}}_{n:n})  \}$.

Here, $\pi_{a}(\cdot; \cdot)=\sum_{i \in \mathcal{I}}\pi_{a,i}(\cdot; \cdot)$. The vector $\hat{\boldsymbol{d}}_{n:n}$ denotes a vector of storage actions in which storages $i<n$ play their NE and storage $n $ cooperates. Similarly, $\hat{\boldsymbol{\tau}}_{n:n}$ denotes a vector of aggregator actions in which it plays the NE with storages $i<n$ and cooperates with storage $n$. The vector $\boldsymbol{\tau}^{D_n}$ represents the action of the aggregator where $\boldsymbol{\tau}^{D_n}_i=\doverline{\boldsymbol{\tau}}_i \; \forall i<n$ (\emph{i.e.} NE),  $\boldsymbol{\tau}^{D_n}_i=\boldsymbol{\tau}_i^D$ if $i=n$ (\emph{i.e.} defection strategy).

Similarly, cooperation is sustained by storage $i$, or equivalently, $w_i^*=\infty $ when $\hat{\boldsymbol{\tau}}_i$ is in the set $\mathcal{A}_{i}^{\mathrm{s}} = \{ \boldsymbol{\tau}_i | \left(1-\delta\right)\pi_i(\boldsymbol{d}_i^D; \boldsymbol{\tau}_i\!) \le\pi_i\left(\boldsymbol{d}_i^*; \boldsymbol{\tau}_i\! \right)\}$. 

When  $\hat{\boldsymbol{\tau}}_i$ is at the boundary of $\mathcal{A}_{i}^{\mathrm{a}}$ ($\mathcal{A}_{i}^{\mathrm{s}}$), the aggregator (storage) is indifferent between cooperating and not cooperating with storage $i$. We assume that when indifferent, both the aggregator and storage cooperate. Both parties cooperate when $\hat{\boldsymbol{\tau}}_i \in \mathcal{A}_{i}$ where $\mathcal{A}_{i}=\mathcal{A}_{i}^{\mathrm{a}} \cap \mathcal{A}_{i}^{\mathrm{s}} $.

\begin{figure}[ht]
  \centering
   \includegraphics[width=0.47\textwidth]{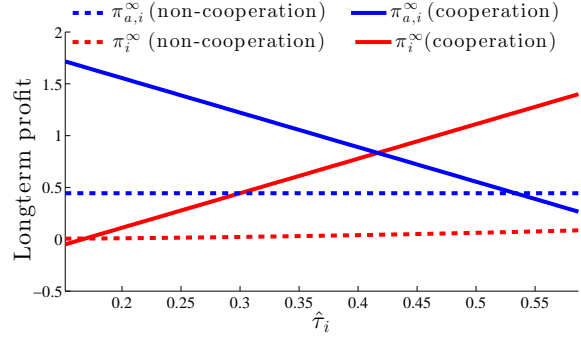}
      \caption{Longterm profits of the aggregator ($\pi_{a,i}^\infty$) and storage $(\pi_i^\infty$) as a function of agreed $\hat{\tau}_i$. For the dashed red line, the storage cheats while the aggregator cooperates. Conversely, for the dashed blue line, the aggregator cheats while the storage cooperates. Both the storage and the aggregator cooperates when $\hat{\tau}_i$ is such that the cooperation profit is greater than the non-cooperation profit for both players.} \label{Figure2}
\end{figure}

\subsubsection{Numerical example} 
Consider a repeated aggregator-storage game with the same aggregator and storages described in the previous examples. We assume that the discount rate is $\delta=0.95$. As it will be shown in Lemma 6, the aggregator and storages agree on the GC solution: $\hat{d}_i =d_i^*=\frac{1}{6}\; \forall i=1,2$. As mentioned in the proof of Lemma 5, we assume that both the aggregator and storages defect from cooperation by maximizing their single game profits so $\tau_i^D = -2 $ and $d_i^D = \frac{\hat{\tau}_i}{2} \; \forall i=1,2$.

The aggregator maximizes its longterm profit by choosing values for $v_i \in \mathbb{R}_+\; \forall\; i=1,2$ such that its longterm profit
 \begin{align*}
  \pi^{\infty}_{a,i} = \frac{\frac{1}{9} - \frac{1}{6}\hat{\tau}_i + 0.95^{v_i} \left( \frac{\hat{\tau}_2}{6}  - \frac{4}{45}\right) }{1-0.95} \quad \forall \; i=1,2
 \end{align*}
 is maximized.
 The aggregator maximizes its longterm profit by choosing to not cooperate (\emph{i.e.} $v_i=0$) if $ \frac{\hat{\tau}_i}{6} - \frac{4}{45}$ is greater than zero. And chooses to cooperate indefinitely (\emph{i.e.}$v_i=\infty$) if $ \frac{\hat{\tau}_i}{6} - \frac{4}{45}$ is less than zero. When $ \frac{\hat{\tau}_i}{6}  - \frac{4}{45} = 0$ the aggregator is indifferent between cooperating and not cooperating. We assume that when indifferent, the aggregator cooperates. In order for the aggregator to have the incentive to cooperate in the longterm, the aggregator must agree to a price schedule such that $\tau_i \lesssim 0.53 \; \forall i=1,2$.
 
 Similarly, each storage maximizes its longterm profit by choosing values for $w_i \in \mathbb{R}_+$ such that 
 \begin{align*}
  \pi^{\infty}_i = \frac{-\frac{1}{36} + \frac{1}{6} \hat{\tau}_i + 0.95^{w_i} \left( \frac{\hat{\tau}_i^2}{80} - \frac{\hat{\tau}_i}{6} +\frac{1}{36}\right)}{1-0.95}
 \end{align*}
 is maximized. Storage $i$ cooperates (\emph{i.e.} chooses $w_i=\infty$) when $\frac{\hat{\tau}_i^2}{80} - \frac{\hat{\tau}_i}{6} +\frac{1}{36} \le 0$ or equivalently when $0.17 \lesssim \hat{\tau}_i \lesssim 13.17$. It follows every storage and the aggregator cooperate when the agreed price schedule satisfies $0.17 \lesssim \hat{\tau}_i \lesssim 0.53 \; \forall i=1,2$. This range of price schedules that foster cooperation between all players can be interpreted as a share in profits between the aggregator and the storages. In the subsection that follows we use Nash's axiomatic bargaining model to predict that the aggregator and the storage will agree on the GC solution and that the profit will be split equally among them when all players are risk neutral. 

\subsubsection*{Proof of Lemma 5}
First we show that for the aggregator to sustain cooperation with storage $i$, $\hat{\boldsymbol{\tau}}_i$ must be in the set $\mathcal{A}_i^{\mathrm{a}}$. The aggregator cooperates with storage $i$ (\emph{i.e.} sends agreed $\hat{\boldsymbol{\tau}}_i$) until time $v_i$, when it cheats (\emph{i.e.} sends defection strategy $\tau_i^D$). Then the repeated game profit for the aggregator $\pi^{\infty}_{a}$  can be expressed as
\begin{align}
&\pi^{\infty}_{a} \!\!= \!\!\sum_{j=1}^n{\!\left\{ \delta^{v_j} \pi_{a}(\boldsymbol{\tau}^{D_j}; \hat{\boldsymbol{d}}_{j:n})+ \!\!\!\!\!\! \sum_{k=v_{j-1}+1}^{v_j-1}{ \!\!\!\!\! \delta^k \pi_{a}(\hat{\boldsymbol{\tau}}_{j:n}; \hat{\boldsymbol{d}}_{j:n}  )  } \right\}}  \label{proof_Lemma_n6_1}
\end{align}
 where $\boldsymbol{\tau}^{D_j}$ denotes the action of the aggregator where $\boldsymbol{\tau}^{D_j}_i=\doverline{\boldsymbol{\tau}}_i \; \forall i<j$ (\emph{i.e.} NE),  $\boldsymbol{\tau}^{D_j}_i=\boldsymbol{\tau}_i^D$ if $i=j$ (\emph{i.e.} defection strategy), and $\boldsymbol{\tau}^{D_j}_i=\hat{\boldsymbol{\tau}}_i \; \forall i>j$ (\emph{i.e.} cooperation strategy), $\hat{\boldsymbol{d}}_{j:n}$ denotes a vector of storage actions in which storages $i<j$ play their NE and storages $i \ge j$ cooperate. Similarly, $\hat{\boldsymbol{\tau}}_{j:n}$ denotes a vector of aggregator actions in which it plays the NE with storages $i<j$ and cooperate with $i\ge j$. It is assumed that the set of storages $\mathcal{I}$ is ordered such that $v_1\le v_2 \le \hdots \le v_n$ and define $v_0\equiv -1$.

We can use the identity $\sum_{k=a+1}^{b-1}{\delta^k}=\frac{\delta^{a+1}-\delta^b}{1-\delta}$ to rewrite \eqref{proof_Lemma_n6_1} as
\begin{align}
\pi^{\infty}_{a} \!\!&=\!\!\! \sum_{j=1}^n{\!\!\{ \delta^{v_j} \pi_{a}(\boldsymbol{\tau}^{D_j}; \hat{\boldsymbol{d}}_{j:n})+{ \frac{\delta^{v_{j\!-
\!1}\!+\!1 }\!\!-\!\delta^{v_{j} }  }{1-\delta} \pi_{a}(\hat{\boldsymbol{\tau}}_{j:n}; \hat{\boldsymbol{d}}_{j:n}  )  } \}}  \nonumber \\
&=\frac{1}{1\!-\!\delta}\bigg( \sum_{j=1}^n{\!\!  \delta^{v_{j}\! }   [ (1-\delta)\pi_{a}(\boldsymbol{\tau}^{D_j}; \hat{\boldsymbol{d}}_{j:n}) -   \pi_{a}(\hat{\boldsymbol{\tau}}_{j:n}; \hat{\boldsymbol{d}}_{j:n}  )    ]    } \nonumber \\
& +  \sum_{j=0}^{n-1} \delta^{v_j\!+\!1} \pi_{a}(\hat{\boldsymbol{\tau}}_{j+1:n}; \hat{\boldsymbol{d}}_{j+1:n}  )   \bigg). \label{Lemma_7_proof_1}
\end{align}

Since $v_0=-1$, $\delta^{v_0+1}=\delta^{-1+1}=1$ and defining $\pi_{a}(\hat{\boldsymbol{\tau}}_{n+1:n} ; \hat{\boldsymbol{d}}_{n+1:n} ) \equiv 0$, equation \eqref{Lemma_7_proof_1} can be rewritten as
\begin{align*}
&\pi^{\infty}_{a} \!\!=\frac{1}{1\!-\!\delta}\left(  \pi_{a}(\hat{\boldsymbol{\tau}} ; \hat{\boldsymbol{d}} )  + \sum_{j=1}^n\!\!  \delta^{v_{j}\! }   \alpha_j  \right)  
\end{align*}
where $\alpha_j=(1\!-\!\delta)\pi_{a}(\boldsymbol{\tau}^{D_j}; \hat{\boldsymbol{d}}_{j:n}) -  \! \pi_{a}(\hat{\boldsymbol{\tau}}_{j:n}; \hat{\boldsymbol{d}}_{j:n}  )  \!+\! \delta \pi_{a}(\hat{\boldsymbol{\tau}}_{j+1:n}; \hat{\boldsymbol{d}}_{j+1:n} )$. Notice that, since $\lim_{v_j \rightarrow \infty}{\delta^{v_j}=0}$ and the aggregator maximizes longterm profits, for the aggregator to choose $v_j=\infty$, $\alpha_j$ must be less than or equal to zero.

Of all possible orderings of the set $\mathcal{I}$ and $d_i^*\ge0\; \forall \; i \in \mathcal{I}$, $\hat{\tau}_j$ it is most constrained when $j=n$, that is, when the least amount of players are in the market\footnote{If the aggregator has defected all previous $n-1$ player, then there are more opportunities for arbitrage in the market (\emph{i.e.} a larger profit to be shared with the storage) and therefore a larger upside when cheating in the single shot game.}. Thus, in order for the aggregator to maintain cooperation with storage $n$, $ \pi_{a}(\hat{\boldsymbol{\tau}}_{n:n}; \hat{\boldsymbol{d}}_{n:n}  )\ge(1\!-\!\delta)\pi_{a}(\boldsymbol{\tau}^{D_n}; \hat{\boldsymbol{d}}_{n:n}) $. The storage maintains cooperation with all storages when the previous equation holds true for any ordering of the set $\mathcal{I}$.

Now we show that for storage $i$ to sustain cooperation with the aggregator, $\hat{\boldsymbol{\tau}}_i$ must be in the set $\mathcal{A}_i^{\mathrm{s}}$. If we assume that storage $i$ cooperates (\emph{i.e.} plays $\hat{\boldsymbol{d}}_i$) until it decides to cheat at time $w_i$ (\emph{i.e.} plays defection strategy $\boldsymbol{d}_i^D$), then the longterm profit of storage $i$, $ \pi^{\infty}_i = \sum_{k=0}^{\infty}{ \delta^k\pi_i(\boldsymbol{d}_i^{(\!k\!)}; \boldsymbol{\tau}_i^{(\!k\!)})  } $ can be expressed as
\begin{align}
\pi^{\infty}_i \!\!&=\!\!\!\! \sum_{k=0}^{w_i-1}{ \!\! \delta^k\pi_i(\hat{\boldsymbol{d}}_i; \hat{\boldsymbol{\tau}}_i)  } \!+\! \delta^{w_i} \pi_i(\boldsymbol{d}_i^D; \hat{\boldsymbol{\tau}}_i\!)\!+\!\!\!\!\!\!\!\sum_{k=w_i+1}^{\infty}{ \!\!\!\!\!\!\delta^k\pi_i(\doverline{\boldsymbol{d}}_i; \doverline{\boldsymbol{\tau}}_i)  } . \label{repeated_game_profit_storage_2} 
\end{align}

Using the identity $\sum_{k=a+1}^{b-1}{\delta^k}=\frac{\delta^{a+1}-\delta^b}{1-\delta}$, equation \eqref{repeated_game_profit_storage_2} can be expressed as
\begin{align*}
\pi^{\infty}_i \!\!&=\frac{1-\delta^{w_i}}{1-\delta}  \pi_i(\hat{\boldsymbol{d}}_i; \hat{\boldsymbol{\tau}}_i)   \!+\! \delta^{w_i} \pi_i(\boldsymbol{d}_i^D; \hat{\boldsymbol{\tau}}_i)\!  \\
&=\! \frac{1}{1\!\!-\!\delta}\!\left(  \!\pi_i\!\left(\!\hat{\boldsymbol{d}}_i; \hat{\boldsymbol{\tau}}_i\!\right)\!  \!+ \! \delta^{w_i}  \!\!\left(\!\left(1\!-\!\delta\right)\!\pi_i\!\left(\!\boldsymbol{d}_i^D; \hat{\boldsymbol{\tau}}_i\!\right)\!-\!\pi_i\!\left(\hat{\boldsymbol{d}}_i; \hat{\boldsymbol{\tau}}_i\!\right)\!\right)
\! \right)
\end{align*}

For storage $i$ to indefinitely sustain cooperation with the aggregator, the term $(1-\delta)\pi_i(\boldsymbol{d}_i^D; \hat{\boldsymbol{\tau}}_i\!)-\pi_i(\hat{\boldsymbol{d}}_i; \hat{\boldsymbol{\tau}}_i)$ must be less or equal to zero. It follows that for both storages and aggregator to sustain cooperation, $\hat{\boldsymbol{\tau}}_i$ must be in both $\mathcal{A}_i^{\mathrm{s}}$  and $\mathcal{A}_i^{\mathrm{a}}$. \QEDB


\subsection{Profit split via Nash Bargaining}  
\label{section:NB}
As shown previously, there are potentially infinitely many ways to split the profit and ensure cooperation. In this section we use John Nash's bargaining model to predict the profit split between the aggregator and the storage. We restrict our analysis to cases where $\mathcal{A}_i \ne \emptyset\; \forall i \in \mathcal{I}$.\footnote{The set $\mathcal{A}_i$ could be empty due to a combination of the following situations: a) the market is too crowded b) the value of money depreciates too rapidly to ensure cooperation or c) the profitability of cheating is too great. }



\subsubsection{The Nash bargaining problem}

In this subsection, we introduce Nash's axiomatic approach to bargaining \cite{Nash-1950}. Denote the set of possible bargaining outcomes of the aggregator and storage $i$ as $\mathcal{S}_i$ and the solution to the bargaining problem as $\xi\left(\mathcal{S}_i\right)$. Under Nash's assumptions a bargaining solution is a single point in a set of possible outcomes that satisfies the following axioms:

\begin{itemize}
\item \emph{Pareto efficiency}: Let $u^{\mathrm{a}}_i(\cdot)$ and  $u^{\mathrm{s}}_i(\cdot)$ be the utility functions of the aggregator and storage $i$ respectively. If $\boldsymbol{a},\;\boldsymbol{b} \in \mathcal{S}_i$, $u^{\mathrm{a}}_i(\boldsymbol{a})>u^{\mathrm{a}}_i(\boldsymbol{b})$, and $u^{\mathrm{s}}_i(\boldsymbol{a})>u^{\mathrm{s}}_i(\boldsymbol{b})$ then $\boldsymbol{b} \ne\xi\left(\mathcal{S}_i\right)$. 
\item \emph{Independence of irrelevant alternatives}: If $\mathcal{S}_i \subseteq \mathcal{T}_i$ and $\xi\left(\mathcal{T}_i\right) \in \mathcal{S}_i$, then $\xi\left(\mathcal{S}_i\right) = \xi\left(\mathcal{T}_i\right)$.
\item \emph{Symmetry}: If $\mathcal{S}_i$ is symmetric ($\exists \; u^{\mathrm{a}}_i(\cdot), \; u^{\mathrm{s}}_i(\cdot)$ such that if $(a,b) \in \mathcal{S}_i$ then $(b,a) \in \mathcal{S}_i$) and $u^{\mathrm{a}}_i(\cdot)$ and $u^{\mathrm{s}}_i(\cdot)$ exhibit this, then the bargaining solution has the form $\boldsymbol{a}=(a,a) = \xi\left(\mathcal{S}_i\right)$ and $u^{\mathrm{s}}_i(\boldsymbol{a})=u^{\mathrm{a}}_i(\boldsymbol{a})$.
\end{itemize}

In this analysis, the aggregator is assumed to independently bargain with each storage. This assumption is reasonable because, by the Pareto efficiency axiom, the agreed storage actions will be the GC solution (\emph{i.e.} $\hat{\boldsymbol{d}}_i=\boldsymbol{d}_i^*)$. Then, the only thing that is left for negotiation is the price schedules sent to each storage. The price schedule sent to a storage does not affect other storages.

\subsubsection*{Lemma 6} 

If both players are risk neutral the aggregator and storage $i$ will agree on a $\hat{\boldsymbol{\tau}}_i$ that equally splits the GC profit and fosters longterm cooperation.

By agreeing to act under the coordination of an aggregator, the storages share some of the profit with the aggregator who is essentially a middle man. However, as seen in figure \ref{Figure2}, as the number of storages increases, splitting the GC profit with an aggregator rather than obtaining the NE profit becomes increasingly lucrative.

It is worth noting that, as shown by \cite{Binmore}, if one of the players is more risk adverse than the other, its share of the profit will decrease. Conversely, if a player is more risk-loving than the other, its share of the profit will increase. 

\subsubsection{Numerical example}

By the Pareto efficiency axiom, $\hat{d}_i = d_i^* =\sfrac{1}{6}\; \forall i=1,2$. Any other choice of $\hat{d}_i$ will yield a lower total profit and could be improved without any player being affected by choosing instead $d_i^*$. Similarly, by Pareto efficiency, $\hat{\tau}_i$ will be one that fosters longterm cooperation $0.17 \lesssim \hat{\tau}_i \lesssim 0.53 \; \forall i=1,2$. 

It follows that the aggregator longterm profit that the aggregator derives from trading with storage $i$ is $\pi_{a,i}^\infty = \sfrac{\left(\frac{1}{9} - \frac{1}{6}\hat{\tau}_i\right)}{\left(1-0.95\right)} $ where $0.17 \lesssim \hat{\tau}_i \lesssim 0.53 \; \forall i=1,2$. Similarly, the longterm profit for storage $i$ is $\pi_i^\infty= \sfrac{\left(-\frac{1}{36} + \frac{1}{6}\hat{\tau}_i\right)}{\left(1-0.95\right)} $ where $0.17 \lesssim \hat{\tau}_i \lesssim 0.53 \; \forall i=1,2$. 

We can express the aggregator's longterm profit derived from trading with storage $i$ as a function of the profit of storage $i$ and normalize it by their joint longterm profit as follows: 
\begin{align*}
\tilde{\pi}_{a,i}^\infty  = 1 - \tilde{\pi}_{i}^\infty \quad \mathrm{where} \; 0 \le \tilde{\pi}_{i}^\infty \le \frac{11}{15}
\end{align*}
where $\tilde{\pi}_{a,i}^\infty$ and $\tilde{\pi}_{i}^\infty$ are the aggregator and storage longterm profit, respectively, normalized by the joint longterm profit $\frac{5}{3}$.

From \cite{Binmore} and assuming that all players are risk-neutral (\emph{i.e.} their utility function is equal to their profit), the solution to the bargaining problem is given by 
\begin{align*}
\mathrm{deal}&= \argmax_{0 \le \tilde{\pi}_{i}^\infty \le \frac{11}{15}} { \left\{ \tilde{\pi}_{i}^\infty \left( 1-\tilde{\pi}_{i}^\infty \right)\right\}} = \sfrac{1}{2}
\end{align*} 
which maps to a price schedule of $\hat{\tau}_i = \sfrac{5}{12}$.

\subsubsection*{Proof of Lemma 6} 
By Pareto efficiency, $\hat{\boldsymbol{d}}_i=\boldsymbol{d}_i^*$ and $\hat{\boldsymbol{\tau}}_i \in \mathcal{A}_i$. Since we assume that the aggregator and storages are risk neutral, their utilities are equal to their profits. The longterm profit of the aggregator from trading with storage $i$ is
\begin{align*}
u_i^{\mathrm{a}}(\hat{\boldsymbol{\tau}}_i)= \sum_{k=0}^\infty{\delta^k }\pi_{a,i}( \hat{\boldsymbol{\tau}}_i;\boldsymbol{d}_i^*\!)= \frac{\pi_{a,i}( \hat{\boldsymbol{\tau}}_i;\boldsymbol{d}_i^*\!)  }{1-\delta} \quad \forall \; \hat{\boldsymbol{\tau}}_i \in  \mathcal{A}_i^{\mathrm{s}}
\end{align*}
and the utility of storage $i$
\begin{align*}
u_i^{\mathrm{s}}(\hat{\boldsymbol{\tau}}_i)=\sum_{k=0}^\infty{\delta^k }\pi_i(\boldsymbol{d}_i^*; \hat{\boldsymbol{\tau}}_i\!)= \frac{\pi_i(\boldsymbol{d}_i^*; \hat{\boldsymbol{\tau}}_i\!)  }{1-\delta}\quad \forall\;\hat{\boldsymbol{\tau}}_i \in  \mathcal{A}_i^{\mathrm{s}}
\end{align*}
is the longterm profit from cooperating with the aggregator. 

Define the set possible agreement outcomes between the aggregator and storage $i$ as
\begin{align*}
 \tilde{\mathcal{S}}_i^{\mathrm{deal}}&=\left\{(\; u_i^{\mathrm{a}}(\hat{\boldsymbol{\tau}}_i),  u_i^{\mathrm{s}}(\hat{\boldsymbol{\tau}}_i)  \;)|\; \hat{\boldsymbol{\tau}}_i \in \mathcal{A}_{i} \right\}  
\end{align*}
and the non-agreement outcomes as $ \tilde{h}_i= (\tilde{h}_i^\mathrm{a},\tilde{h}_i^\mathrm{s}) = (0,0)$. The set of possible outcomes is then $\tilde{\mathcal{S}_i} = \tilde{\mathcal{S}}_i^{\mathrm{deal}} \cup h_i$. Since $\tilde{h}_i \notin \tilde{\mathcal{S}}_i^{\mathrm{deal}}$, $\tilde{\mathcal{S}_i}$ is not necessarily convex. We define the possible bargaining outcomes as the convex hull of $\tilde{\mathcal{S}_i}$, (\emph{i.e.} $\mathcal{S}_i=\mathrm{conv}(\tilde{\mathcal{S}_i})$).

From the Pareto efficiency axiom, we know that the deal will lie on $\tilde{\mathcal{S}}_i^{\mathrm{deal}}$. From the symmetry axiom, we know that $u_i^{\mathrm{a}}(\xi\left(\mathcal{S}_i\right))=u_i^{\mathrm{s}}(\xi\left(\mathcal{S}_i\right))$. Hence, the aggregator and the storage choose a $\hat{\boldsymbol{\tau}}_i$ that equally splits the profit. 

A way of defining the solution function $\xi\left(\mathcal{S}_i\right)$ is given by Binmore in \cite{Binmore}. Binmore shows that $\xi\left(\mathcal{S}_i\right)= \argmax_{( u_i^{\mathrm{a}}(\hat{\boldsymbol{\tau}}_i),  u_i^{\mathrm{s}}(\hat{\boldsymbol{\tau}}_i) ) \in \mathcal{S}_i \ge \tilde{h}_i} {( u_i^{\mathrm{a}}(\hat{\boldsymbol{\tau}}_i) - h_i^\mathrm{a})( u_i^{\mathrm{s}}(\hat{\boldsymbol{\tau}}_i) - h_i^\mathrm{s})} $. From the Pareto efficiency axiom we know that $\xi\left(\mathcal{S}_i\right) \in  \tilde{\mathcal{S}}_i^{\mathrm{deal}}$. Substituting $(h_i^\mathrm{a},h_i^\mathrm{s}) = (0,0)$ and from the independence of irrelevant alternatives axiom, we arrive at 
\begin{align}
\xi\left(\mathcal{S}_i\right)&= \argmax_{( u_i^{\mathrm{a}}(\hat{\boldsymbol{\tau}}_i),  u_i^{\mathrm{s}}(\hat{\boldsymbol{\tau}}_i)) \in \tilde{\mathcal{S}}_i^{\mathrm{deal}} \ge (0,0)} {u_i^{\mathrm{a}}(\hat{\boldsymbol{\tau}}_i)  u_i^{\mathrm{s}}(\hat{\boldsymbol{\tau}}_i)} \label{eq:binmore}
\end{align}

Notice that we can write the longterm utility of the aggregator from cooperating with storage $i$ as a function of the storage's utility

\begin{align*}
u_i^{\mathrm{a}}(\hat{\boldsymbol{\tau}}_i) \!=\! \underbrace{\frac{1}{1-\delta}\sum_{t \in \mathcal{T}} \! \left\{ \!\!-p^{[t]}\!\!\left(\boldsymbol{d}^{[t]*}\right) \cdot d_i^{[t]*}  \!-\! \frac{\epsilon_i}{2}d_i^{[t]*2}\!\right\} }_{\pi_\mathrm{total}=\text{Total profit to be shared}}
\!-u_i^{\mathrm{s}}(\hat{\boldsymbol{\tau}}_i).
\end{align*}
Normalizing by the total profit to be shared, $\pi_\mathrm{total}$, the longterm utility of the aggregator from cooperating with storage $i$ can be written as: $\tilde{u}_i^{\mathrm{a}}(\hat{\boldsymbol{\tau}}_i) \!=1-\tilde{u}_i^{\mathrm{s}}(\hat{\boldsymbol{\tau}}_i)$, where $\tilde{u}_i^{\mathrm{a}}(\hat{\boldsymbol{\tau}}_i) =\sfrac{u_i^{\mathrm{a}}(\hat{\boldsymbol{\tau}}_i)}{\pi_\mathrm{total}} $ and $\tilde{u}_i^{\mathrm{s}}(\hat{\boldsymbol{\tau}}_i) =\sfrac{u_i^{\mathrm{s}}(\hat{\boldsymbol{\tau}}_i)}{\pi_\mathrm{total}} $. 

We can rewrite \eqref{eq:binmore} as
\begin{align*}
\xi\left(\mathcal{S}_i\right)&= \argmax_{( \tilde{u}_i^{\mathrm{s}}(\hat{\boldsymbol{\tau}}_i),  1-\tilde{u}_i^{\mathrm{s}}(\hat{\boldsymbol{\tau}}_i)) \in \tilde{\mathcal{S}}_i^{\mathrm{deal}} \ge (0,0)} {  \tilde{u}_i^{\mathrm{s}}(\hat{\boldsymbol{\tau}}_i)}\left(  1-\tilde{u}_i^{\mathrm{s}}(\hat{\boldsymbol{\tau}}_i)\right).
\end{align*}
whose solution is $\tilde{u}_i^{\mathrm{s}}(\hat{\boldsymbol{\tau}}_i)=\frac{1}{2}$ and thus $\tilde{u}_i^{\mathrm{a}}(\hat{\boldsymbol{\tau}}_i)=\frac{1}{2}$. \QEDB

\section{Conclusion}
\label{section: Conclusion}

We studied the profit of a group of energy storages under competition and cooperation. We showed that without cooperation, the aggregate profit of the storages approaches zero as the number of storages grows. We presented two approaches to foster cooperation. In the first approach, storages are exposed to artificial cost functions, and their self-interested strategy maximizes the aggregate profit. In the second approach, the aggregate profit is maximized with the help of an aggregator. The interaction of the aggregator and storages is modeled as a simultaneous move game whose Nash equilibrium is undesirable. We derive the conditions (\emph{i.e.} profit split between aggregator and storages) that ensure aggregator-storage cooperation. Finally, we use Nash's axiomatic approach to bargaining to predict that risk-neutral players will equally split the available profit.

\bibliographystyle{IEEEtran}
\bibliography{bibliography}

\end{document}